# Some remarks about q- Chebyshev polynomials and q- Catalan numbers and related results

Johann Cigler

Fakultät für Mathematik, Universität Wien
Uni Wien Rossau, Oskar-Morgenstern-Platz 1, 1090 Wien
johann.cigler@univie.ac.at
http://homepage.univie.ac.at/johann.cigler/

**Abstract**

The moments of the Lucas polynomials and of the Chebyshev polynomials of the first kind are (multiples of) central binomial coefficients and the moments of the Fibonacci polynomials and of the Chebyshev polynomials of the second kind are Catalan numbers. In this survey paper we present some generalizations of these results together with various $q-$ analogues.

**0. Introduction**

The moments of the Fibonacci polynomials and of the Chebyshev polynomials of the second kind are (multiples of) the Catalan numbers $C_n = \frac{1}{n+1}\binom{2n}{n}$ and the moments of the Lucas polynomials or equivalently of the Chebyshev polynomials of the first kind are (multiples of) the central binomial coefficients $B_n = \binom{2n}{n}$. We show first how these facts generalize to various $q-$ analogues. The most natural $q-$ analogues of the monic Chebyshev polynomials (cf. [11] or [12]) are orthogonal polynomials and their moments are multiples of 
$c_n(q) = \frac{1}{[n+1]}\begin{bmatrix}2n\\n\end{bmatrix}$ and $b_n(q) = \begin{bmatrix}2n\\n\end{bmatrix}$ respectively. The same fact holds for a curious class (cf. [8]) of non-orthogonal $q-$Fibonacci and $q-$Lucas polynomials, whereas the moments $C_n(q)$ of the Carlitz $q-$ Fibonacci polynomials - which are orthogonal - do not have explicit expressions. But the generating function $C_q(u) = \sum_{n\geq 0} C_n(q) u^n$ has a simple representation as $C_q(u) = \frac{E(-qu)}{E(-u)}$ for some $q-$ analogue $E(u)$ of the exponential series. Some of these results can be extended to $q-$ analogues of generalized Fibonacci polynomials $f_n^{(m)}(x)$ whose moments are multiples of generalized Catalan numbers $C_n^{(m)} = \frac{1}{(m-1)n+1}\binom{mn}{n}$ and to $q-$ analogues of generalized Lucas polynomials whose moments are related to generalized central binomial coefficients $\binom{mn}{n}$.

For the convenience of the reader I first recall some well-known background material about Fibonacci and Lucas polynomials and their generalizations $f_n^{(m)}(x)$ and $l_n^{(m)}(x)$ and state some basic facts about $q-$ identities.

## 1. Background material

**1.1.** Let the special Fibonacci polynomials $f_n(x)$ be defined by

$$f_n(x) = \sum_{k=0}^{\lfloor \frac{n}{2} \rfloor} \binom{n-k}{k} (-1)^k x^{n-2k}. \tag{1.1}$$

They satisfy the recurrence relation

$$f_n(x) = x f_{n-1}(x) - f_{n-2}(x) \tag{1.2}$$

with initial values $f_0(x) = 1$ and $f_1(x) = x$ and are orthogonal with respect to the linear functional $\Lambda_f$ defined by $\Lambda_f(f_n(x)) = [n=0]$. More precisely we have

$$\Lambda_f(f_n(x) f_m(x)) = [n = m]. \tag{1.3}$$

The moments $\Lambda_f(x^n)$ can easily be deduced from the formula

$$x^n = \sum_{k=0}^{\lfloor \frac{n}{2} \rfloor} \left( \binom{n}{k} - \binom{n}{k-1} \right) f_{n-2k}(x). \tag{1.4}$$

This gives $\Lambda_f(x^{2n+1}) = 0$ and

$$\Lambda_f(x^{2n}) = C_n, \tag{1.5}$$

where $C_n = \binom{2n}{n} - \binom{2n}{n-1} = \frac{1}{n+1}\binom{2n}{n}$ is a Catalan number.

The moment generating function is

$$\sum_{n \geq 0} \Lambda_f(x^{2n}) u^n = C(u) = \frac{1 - \sqrt{1-4u}}{2u} \tag{1.6}$$

and satisfies

$$C(u) = 1 + u C(u)^2. \tag{1.7}$$

Equivalently we have

$$\frac{1}{1+u} C\left( \frac{u}{(1+u)^2} \right) = 1. \tag{1.8}$$

If we set $\alpha = \dfrac{x + \sqrt{x^2 - 4}}{2}$ and $\beta = \dfrac{x - \sqrt{x^2 - 4}}{2}$, then $\alpha^n = x\alpha^{n-1} - \alpha^{n-2}$ and $\beta^n = x\beta^{n-1} - \beta^{n-2}$. Note that $\alpha + \beta = x$ and $\alpha\beta = 1$.

Since $\dfrac{\alpha^{n+1} - \beta^{n+1}}{\alpha - \beta}$ satisfies recursion (1.2) and the initial values we get the Binet formulae

$$f_n(x) = \frac{\alpha^{n+1} - \beta^{n+1}}{\alpha - \beta}. \tag{1.9}$$

**1.2.** Let us also consider a variant $l_n(x)$ of the Lucas polynomials defined by

$$l_n(x) = \sum_{k=0}^{\left\lfloor \frac{n}{2} \right\rfloor} \binom{n-k}{k} \frac{n}{n-k} (-1)^k x^{n-2k} \qquad (1.10)$$

$$l_0(x) = 1.$$

The polynomials $l_n(x)$ satisfy the recurrence relation

$$l_n(x) = x l_{n-1}(x) - \lambda_{n-2} l_{n-2}(x) \qquad (1.11)$$

with $\lambda_0 = 2$ and $\lambda_n = 1$ for $n > 0$.

We have $l_n(x) = f_n(x) - f_{n-2}(x)$ for $n \geq 2$ and $l_n(x) = f_n(x)$ for $n = 0, 1$.

This implies for $n > 0$ the Binet formulae

$$l_n(x) = \alpha^n + \beta^n. \qquad (1.12)$$

The polynomials $l_n(x)$ are orthogonal with respect to the linear functional $\Lambda_l$ defined by $\Lambda_l(l_n(x)) = [n = 0]$. More precisely

$$\Lambda_l(l_n(x) l_m(x)) = 2[n = m] \qquad (1.13)$$

for $n > 0$ and $\Lambda_l(l_0(x)^2) = 1$.

From the representation

$$x^n = \sum_{k=0}^{\left\lfloor \frac{n}{2} \right\rfloor} \binom{n}{k} l_{n-2k}(x) \qquad (1.14)$$

we deduce that the moments $\Lambda_l(x^n)$ are $\Lambda_l(x^{2n+1}) = 0$ and

$$\Lambda_l(x^{2n}) = B_n = \binom{2n}{n}, \qquad (1.15)$$

where $B_n = \binom{2n}{n}$ is a central binomial coefficient.

The moment generating function is

$$\sum_{n \geq 0} \Lambda_l(x^{2n}) u^n = B(u) = \frac{1}{\sqrt{1-4u}}. \qquad (1.16)$$

Equivalently

$$\sum_{n \geq 0} B_n \frac{u^n}{(1+u)^{2n+1}} = \frac{1}{1-u}. \qquad (1.17)$$

**1.3.** Let us note some generalizations of the foregoing situation.

Let $m \geq 1$ and let

$$f_n^{(m)}(x) = \sum_{k=0}^{\left\lfloor \frac{n}{m} \right\rfloor} \binom{n-(m-1)k}{k} (-1)^k x^{n-mk}. \qquad (1.18)$$

These polynomials satisfy

$$f_n^{(m)}(x) = x f_{n-1}^{(m)}(x) - f_{n-m}^{(m)}(x) \qquad (1.19)$$

for $n > 0$ with initial values $f_0^{(m)}(x) = 1$ and $f_{-n}^{(m)}(x) = 0$ for $0 < n < m$. This implies $f_n^{(m)}(x) = x^n$ for $0 \leq n < m$.

We have

$$x^n = \sum_{k=0}^{\left\lfloor \frac{n}{m} \right\rfloor} \left( \binom{n}{k} - (m-1)\binom{n}{k-1} \right) f_{n-mk}^{(m)}(x) \qquad (1.20)$$

and therefore the linear functional $\Lambda_{f^{(m)}}$ defined by $\Lambda_{f^{(m)}}\left(f_n^{(m)}(x)\right) = [n=0]$ gives as moments the $m-$ Catalan numbers, which are also called Fuss-Catalan numbers,

$$\Lambda_{f^{(m)}}\left(x^{nm}\right) = \binom{nm}{n} - (m-1)\binom{nm}{n-1} = \frac{1}{(m-1)n+1}\binom{nm}{n} = C_n^{(m)}. \qquad (1.21)$$

For some special cases cf. [14], OEIS A000108, A001764, A002293, A009294.

The generating function satisfies $\Phi_m(u) = \sum_{n \geq 0} C_n^{(m)} u^n = 1 + u \Phi_m(u)^m$.

The most natural generalization of the Lucas polynomials are the polynomials $l_n^{(m)}(x)$ which satisfy $l_n^{(m)}(x) = x l_{n-1}^{(m)}(x) - \lambda_{n-m} l_{n-m}^{(m)}(x)$ with initial values $l_n^{(m)}(x) = x^n$ for $0 \leq n < m$ and $\lambda_0 = m$ and $\lambda_n = 1$ for $n > 0$.

This gives for $n > 0$

$$l_n^{(m)}(x) = \sum_{k=0}^{\left\lfloor \frac{n}{m} \right\rfloor} \binom{n-(m-1)k}{k} \frac{n}{(n-(m-1)k)} (-1)^k x^{n-mk}. \qquad (1.22)$$

It is easily verified that

$$l_n^{(m)}(x) = f_n^{(m)}(x) - (m-1) f_{n-m}^{(m)}(x) \qquad (1.23)$$

for $n > 0$.

Here we have

$$x^n = \sum_{k=0}^{\left\lfloor \frac{n}{m} \right\rfloor} \binom{n}{k} l_{n-mk}^{(m)}(x) \qquad (1.24)$$

and thus for the linear functional $\Lambda_{l^{(m)}}$ defined by $\Lambda_{l^{(m)}}\left(l_n^{(m)}(x)\right) = [n=0]$

$$\Lambda_{l^{(m)}}\left(x^{nm}\right) = \binom{nm}{n}. \qquad (1.25)$$

The generating function $\Psi_m(u) = \sum_{n\geq 0} \binom{mn}{n} u^n$ satisfies

$$\Psi_m(u) = \frac{1}{1 - mu\Phi_m(u)^{m-1}}.$$

It seems that these polynomials do not have interesting $q-$ analogues. Therefore we consider another generalization $L_n^{(m)}(x)$ of the Lucas polynomials which satisfy
$L_n^{(m)}(x) = xL_{n-1}^{(m)}(x) - \lambda_{n-m} L_{n-m}^{(m)}(x)$ with initial values $L_n^{(m)}(x) = x^n$ for $0 \leq n < m$ and $\lambda_0 = 2$ and $\lambda_n = 1$ for $n > 0$.

Here we get for $n > 0$

$$L_n^{(m)}(x) = \sum_{k=0}^{\left\lfloor \frac{n}{m} \right\rfloor} (-1)^k \binom{n-(m-1)k}{k} \frac{n-(m-2)k}{n-(m-1)k} x^{n-mk}, \qquad (1.26)$$

which implies that

$$L_n^{(m)}(x) = f_n^{(m)}(x) - f_{n-m}^{(m)}(x) \qquad (1.27)$$

for $n > 0$.

For this is trivially true for $1 \leq n \leq m-1$ since in this case $L_n^{(m)}(x) = x^n = x^n - 0$. It is also true for $n = m$, for in this case we have
$L_m^{(m)}(x) = xL_{m-1}^{(m)}(x) - 2L_0^{(m)}(x) = x^m - 2 = (x^m - 1) - 1 = f_m^{(m)}(x) - f_0^{(m)}(x)$.

For $n > m$ both sides satisfy the same recurrence relation.

If we define the linear functional $\Lambda_{L^{(m)}}$ by $\Lambda_{L^{(m)}}\left(L_n^{(m)}(x)\right) = [n=0]$ we get in this case

$$x^n = \sum_{k=0}^{\left\lfloor \frac{n}{m} \right\rfloor} \left(\binom{n}{k} - (m-2)\sum_{j=0}^{k-1} \binom{n}{j}\right) L_{n-mk}^{(m)} \qquad (1.28)$$

and thus

$$\Lambda_{L^{(m)}}(x^{nm}) = \binom{nm}{n} - (m-2)\sum_{j=0}^{n-1}\binom{mn}{j}.  \tag{1.29}$$

**1.4.** Let us now state some well-known notations and results for $q-$ identities which will be needed later (cf. e.g. [5]).

We always assume that $0 < q < 1$. Let $(x;q)_n = (1-x)(1-qx)\cdots(1-q^{n-1}x)$ and $(x;q)_\infty = \prod_{j\geq 0}(1-q^j x)$.

Let $[n] = [n]_q = \dfrac{1-q^n}{1-q}$ and $\begin{bmatrix}n\\k\end{bmatrix} = \begin{bmatrix}n\\k\end{bmatrix}_q = \dfrac{(q;q)_n}{(q;q)_k (q;q)_{n-k}}$ for $0 \leq k \leq n$.

The $q-$ binomial coefficients satisfy

$$\begin{bmatrix}n\\k\end{bmatrix} = q^k \begin{bmatrix}n-1\\k\end{bmatrix} + \begin{bmatrix}n-1\\k-1\end{bmatrix} \quad \text{and} \quad \begin{bmatrix}n\\k\end{bmatrix} = \begin{bmatrix}n-1\\k\end{bmatrix} + q^{n-k}\begin{bmatrix}n-1\\k-1\end{bmatrix}.$$

Let $\varepsilon$ be the operator $\varepsilon f(x) = f(qx)$.

Let $D_q$ be the $q-$ differentiation operator defined by $D_q f(x) = \dfrac{f(x)-f(qx)}{x-qx}$.

Then $\varepsilon = 1 + (q-1)xD_q$ since $(1+(q-1)xD_q)x^n = x^n + (q^n-1)x^n = q^n x^n = \varepsilon x^n$.

A simple $q-$ analogue of the binomial theorem is the fact that the so called Rogers-Szegö polynomials $r_n(x,s) = \sum_{k=0}^n \begin{bmatrix}n\\k\end{bmatrix} x^k s^{n-k}$ can be represented as $r_n(x,s) = (x+s\varepsilon)^n 1$. This follows by induction because

$$(x+s\varepsilon)\sum_k \begin{bmatrix}n\\k\end{bmatrix} x^k s^{n-k} = \sum_k \begin{bmatrix}n\\k\end{bmatrix} x^{k+1} s^{n-k} + \sum_k \begin{bmatrix}n\\k\end{bmatrix} q^k x^k s^{n+1-k} = \sum_k \left(\begin{bmatrix}n\\k-1\end{bmatrix} + q^k \begin{bmatrix}n\\k\end{bmatrix}\right) x^k s^{n+1-k}$$

$$= \sum_k \begin{bmatrix}n+1\\k\end{bmatrix} x^k s^{n+1-k}.$$

Note also that $D_q r_n(x,s) = [n]r_{n-1}(x,s)$ because

$$D_q r_n(x,s) = \sum_k \begin{bmatrix}n\\k\end{bmatrix} D_q x^k s^{n-k} = \sum_k \begin{bmatrix}n\\k\end{bmatrix} [k] x^{k-1} s^{n-k} = [n]\sum_k \begin{bmatrix}n-1\\k-1\end{bmatrix} x^{k-1} s^{n-k} = [n] r_{n-1}(x,s)$$

Therefore the Rogers-Szegö polynomials satisfy the recursion

$$r_{n+1}(x,s) = (x+s+(q-1)sxD_q)r_n(x,s) = (x+s)r_n(x,s) + (q^n-1)xsr_{n-1}(x,s).  \tag{1.30}$$

We shall also need the following $q$-hypergeometric version of the binomial theorem (cf. e.g. [3]).

$$\sum_{k\geq 0}\frac{(a;q)_k}{(q;q)_k}x^k = \frac{(ax;q)_\infty}{(x;q)_\infty} \qquad (1.31)$$

for $|x|<1$.

As special cases we will need the following $q-$ analogues of the exponential series

$$\sum_{n\geq 0}\frac{x^n}{(q;q)_n} = \frac{1}{(x;q)_\infty} \qquad (1.32)$$

and

$$\sum_{n\geq 0}\frac{q^{\binom{n}{2}}(-x)^n}{(q;q)_n} = (x;q)_\infty \qquad (1.33)$$

and the simplest binomial theorems

$$\sum_{k=0}^{n}(-1)^k\begin{bmatrix}n\\k\end{bmatrix}q^{\binom{k}{2}}x^k = (x;q)_n \qquad (1.34)$$

and

$$\sum_{k\geq 0}\begin{bmatrix}n+k-1\\k\end{bmatrix}x^k = \frac{1}{(x;q)_n}. \qquad (1.35)$$

## 1.5. Notes

The polynomials $f_n(x)$ are the special case $f_n(x) = Fib_{n+1}(x,-1)$ of the bivariate Fibonacci polynomials $Fib_n(x,s) = \sum_{k=0}^{\lfloor\frac{n-1}{2}\rfloor}\binom{n-k-1}{k}s^k x^{n-1-2k}$ which satisfy
$Fib_n(x,s) = xFib_{n-1}(x,s) + sFib_{n-2}(x,s)$ with initial values $Fib_0(x,s) = 0$ and $Fib_1(x,s) = 1$.

For $n>0$ the polynomials $l_n(x)$ are the special case $l_n(x) = Luc_n(x,-1)$ of the bivariate Lucas polynomials which satisfy $Luc_n(x,s) = \sum_{k=0}^{\lfloor\frac{n}{2}\rfloor}\frac{n}{n-k}\binom{n-k}{k}s^k x^{n-2k}$ for $n>0$ and $Luc_0(x,s) = 2$.

There is also a close connection with the Chebyshev polynomials in the usual form.
The Chebyshev polynomials of the first kind $T_n(x)$ satisfy $T_n(x) = 2xT_{n-1}(x) - T_{n-2}(x)$ with initial values $T_0(x) = 1$ and $T_1(x) = x$ and $T_n(x) = 2^{n-1}Luc_n\left(x,-\frac{1}{4}\right) = \frac{Luc_n(2x,-1)}{2}$.

The monic Chebyshev polynomials of the first kind

$$t_n(x) = \sum_{k=0}^{\left\lfloor \frac{n}{2} \right\rfloor} \frac{n}{n-k} \binom{n-k}{k} \left(-\frac{1}{4}\right)^k x^{n-2k} \tag{1.36}$$

for $n > 0$ and $t_0(x) = 1$ satisfy $t_n(x) = x t_{n-1}(x) - \lambda_{n-2} t_{n-2}(x)$ with $\lambda_0 = \frac{1}{2}$ and $\lambda_n = \frac{1}{4}$.

The Chebyshev polynomials of the second kind $U_n(x)$ satisfy $U_n(x) = 2x U_{n-1}(x) - U_{n-2}(x)$ with initial values $U_0(x) = 1$ and $U_1(x) = 2x$ and $U_n(x) = Fib_{n+1}(2x, -1) = 2^n Fib_{n+1}\left(x, -\frac{1}{4}\right)$.

The monic Chebyshev polynomials of the second kind

$$u_n(x) = \sum_{k=0}^{\left\lfloor \frac{n}{2} \right\rfloor} \binom{n-k}{k} \left(-\frac{1}{4}\right)^k x^{n-2k} \tag{1.37}$$

satisfy $u_n(x) = x u_{n-1}(x) - \frac{1}{4} u_{n-2}(x)$.

Identities (1.4), (1.14), (1.20), (1.24) and (1.28) are special cases of the following situation:
Let $p_n(x)$ be polynomials satisfying

$$p_n(x) = x p_{n-1}(x) - \lambda_n p_{n-m}(x) \tag{1.38}$$

for some integer $m \geq 1$ and initial values $p_n(x) = x^n$ for $0 \leq n < m$.

Then there are uniquely determined coefficients $c(n, k)$ such that

$$x^n = \sum_{k=0}^{n} c(n, k) p_k(x). \tag{1.39}$$

This implies

$$\sum_k c(n,k) p_k(x) = x \sum_k c(n-1, k) p_k(x) = \sum_k c(n-1, k) \left( p_{k+1}(x) + \lambda_{k+m-1} p_{k-m+1}(x) \right)$$
$$= \sum_k p_k(x) \left( c(n-1, k-1) + \lambda_k c(n-1, k+m-1) \right)$$

Therefore we get

$$c(n, k) = c(n-1, k-1) + \lambda_k c(n-1, k+m-1) \tag{1.40}$$

with initial values $c(0, k) = [k = 0]$ and boundary values $c(n, -1) = 0$.

If we apply the linear functional $\Lambda$ we get

$$\Lambda(x^n) = c(n, 0). \tag{1.41}$$

Let us recall a well-known combinatorial interpretation of (1.40):

The number $c(n,k)$ is the weight of all lattice paths in $\mathbb{N}^2$ which initial point $(0,0)$ and endpoint $(n,k)$, where each step is either an up-step $(j,\ell) \to (j+1,\ell+1)$ or a down-step $(j,\ell+m-1) \to (j+1,\ell)$. The weight of an up-step is 1 and the weight of a down-step with endpoint $(*,k)$ is $\lambda_k$. The weight of a path is the product of the weights of its steps. The weight of a set of paths is the sum of their weights. The trivial path $(0,0) \to (0,0)$ has by definition weight 1.

It is then clear that $c(mn+k+i,k) = 0$ for $0 \le i < m$. Therefore (1.39) can also be written in the form

$$x^n = \sum_{k=0}^{n} c(n, n-km) p_{n-km}(x). \tag{1.42}$$

Let $\Phi(u) = \sum_n c(mn,0) u^n$ be the generating function of the moments. Each nontrivial path has a unique decomposition of the following form: For each $i$ with $1 \le i \le m-1$ there is an up-step from height $i-1$ to height $i$ followed by a maximal path which ends at height $i$ and never falls below height $i$, and finally a down-step which ends on height 0 followed by a non-negative path which ends on height 0.

Therefore we get

$$\Phi(u) = 1 + \lambda_0 u \Phi(u) \Phi_1(u) \Phi_2(u) \cdots \Phi_{m-1}(u) \tag{1.43}$$

where $\Phi_i(u)$ denotes the generating function of the moments if in (1.38) $\lambda_n$ is replaced by $\lambda_{n+i}$.

For $m=2$ we get

$$c(2n,0) = \lambda_0 \sum_{k=0}^{n-1} b(2k,0) c(2n-2-2k,0), \tag{1.44}$$

where $b(n,m)$ is the corresponding weight when $\lambda_k$ is replaced by $\lambda_{k+1}$.

Let $\Phi(u) = \sum_n c(2n,0) u^n$ and $\Psi(u) = \sum_n b(2n,0) u^n$.

Then (1.44) is equivalent with

$$\Phi(u) = 1 + \lambda_0 u \Phi(u) \Psi(u). \tag{1.45}$$

For the polynomials $f_n(x)$ we have $\lambda_n = 1$ and therefore (1.45) reduces to (1.7).

Identity (1.4) is equivalent with $c(2n+k,k) = \binom{2n+k}{n} - \binom{2n+k}{n-1}$ and $c(2n+k+1,k) = 0$.

The corresponding matrix $(c(n,k))$ is known as Catalan triangle (cf. [14], OEIS, A053121).

For the polynomials $l_n(x)$ we have $\lambda_0 = 2$ and $\lambda_n = 1$. Therefore $\Psi(u) = C(u)$ and

$\Phi(u) = 1 + 2u\Phi(u)C(u)$, which gives $\Phi(u) = \dfrac{1}{1-2uC(u)} = \dfrac{1}{\sqrt{1-4u}}$.

Formula (1.24) is easily proved by induction. It is trivially true for $n < m$. For $n = m$ we have $l_m^{(m)}(x) = x^m - m$ and (1.24) reduces to $\sum_{k=0}^{1} \binom{m}{k} l_{m-mk}^{(m)}(x) = l_m^{(m)}(x) + m l_0^{(m)}(x) = x^m$.

If it is true for $n-1 \geq m$ then

$$\sum_{k=0}^{\lfloor n/m \rfloor} \binom{n}{k} l_{n-mk}^{(m)}(x) = \sum_{k=0}^{\lfloor n/m \rfloor} \binom{n-1}{k} l_{n-mk}^{(m)}(x) + \sum_{k=0}^{\lfloor n/m \rfloor} \binom{n-1}{k-1} l_{n-mk}^{(m)}(x) = \sum_{k=0}^{\lfloor n/m \rfloor} \binom{n-1}{k} \left( l_{n-mk}^{(m)}(x) + l_{n-m-mk}^{(m)}(x) \right)$$

$$= \sum \binom{n-1}{k} x l_{n-1-mk}^{(m)}(x) = x^n.$$

From (1.23) we get (1.20) and from (1.27) we deduce (1.28).

## 2. The simplest q-analogues

The simplest $q-$ analogues of $f_n^{(m)}(x) = \sum_{k=0}^{\lfloor n/m \rfloor} \binom{n-(m-1)k}{k}(-1)^k x^{n-mk}$ are the polynomials

$$f_n^{(m)}(x,q) = \sum_{k=0}^{\lfloor n/m \rfloor} q^{m\binom{k}{2}} \begin{bmatrix} n-(m-1)k \\ k \end{bmatrix} (-1)^k x^{n-mk}, \tag{2.1}$$

which satisfy

$$f_n^{(m)}(x,q) = x f_{n-1}^{(m)}(x,q) - q^{n-m} f_{n-m}^{(m)}(x,q) \tag{2.2}$$

with initial values $f_0^{(m)}(x,q) = 1$ and $f_{-n}^{(m)}(x,q) = 0$ for $0 < n < m$.

For $m = 2$ these are orthogonal polynomials which are closely related to Carlitz's $q-$ Fibonacci polynomials.

What can be said about the moments of these polynomials?

Denote the moments by $C_n^{(m)}(q)$. The generating function $C_q^{(m)}(u) = \sum_{n \geq 0} C_n^{(m)}(q) u^n$ satisfies

$$C_q^{(m)}(u) = 1 + u C_q^{(m)}(u) C_q^{(m)}(qu) \cdots C_q^{(m)}(q^{m-1}u). \tag{2.3}$$

These $q-$ Catalan numbers have no simple closed formula, but their generating function can be represented (cf. [9]) as

$$C_q^{(m)}(u) = \frac{E^{(m)}(-qu)}{E^{(m)}(-u)} \tag{2.4}$$

with

$$E^{(m)}(u) = \sum_{n\geq 0} \frac{q^{m\binom{n}{2}}}{(q;q)_n} u^n. \tag{2.5}$$

For $E^{(m)}(u) - E^{(m)}(qu) = \sum_{n\geq 0} \frac{q^{m\binom{n}{2}}}{(q;q)_n}(1-q^n)u^n = u\sum_{n\geq 0} \frac{q^{m\binom{n}{2}}}{(q;q)_{n-1}}u^{n-1} = uE^{(m)}(q^m u)$ implies (2.3).

For $m=2$ the generating function $C_q(u) = \sum_{n\geq 0} C_n(q)u^n$ satisfies

$$C_q(u) = 1 + uC_q(u)C_q(qu). \tag{2.6}$$

Comparing coefficients this gives

$C_n(q) = \sum_{k=0}^{n-1} q^k C_k(q) C_{n-1-k}(q)$ with $C_0(q)=1$. Further properties can be found in [13].

From (2.1) we see that

$0 = \Lambda_{f^{(m)},q}\left(f_{nm}^{(m)}(x,q)\right) = \sum_{k=0}^{n} q^{m\binom{k}{2}} \begin{bmatrix} nm-(m-1)k \\ k \end{bmatrix} (-1)^k C_{n-k}^{(m)}(q)$. This can be used to compute $C_n^{(m)}(q)$.

In [6] a shorter algorithm to compute these numbers has been given. We have

$$\sum_{k=0}^{n}(-1)^k q^{2\binom{k}{2}} \begin{bmatrix} n+1-k \\ k \end{bmatrix} C_{n-k}(q) = 0 \tag{2.7}$$

and more generally

$$\sum_{k=0}^{n}(-1)^k q^{m\binom{k}{2}} \begin{bmatrix} 1+(m-1)(n-k) \\ k \end{bmatrix} C_{n-k}^{(m)}(q) = 0. \tag{2.8}$$

Identity (2.7) follows immediately from orthogonality because

$$\sum_{k=0}^{n}(-1)^k q^{2\binom{k}{2}} \begin{bmatrix} n+1-k \\ k \end{bmatrix} C_{n-k}(q) = \Lambda_{f^{(m)},q}\left(x^{n-1} f_{n+1}^{(2)}(x,q)\right) = 0.$$

Identity (2.8) is equivalent with $\Lambda_{f^{(m)},q}\left(x^{n-1} f^{(m)}_{(m-1)n+1}(x,q)\right) = 0.$

Let $r(n,k) = \Lambda_{f^{(m)},q}\left(x^k f^{(m)}_{mn-k}(x,q)\right).$

By (2.2) $r(n,k)$ is a linear combination of $r(n,k-1)$ and $r(n-1,k-1)$.

Therefore $r(n,n-1)$ is a linear combination of $r(n,n-2)$ and $r(n-1,n-2)$ and therefore a linear combination of $r(n,n-3)$, $r(n-1,n-3)$ and $r(n-2,n-3)$ and thus a linear combination of $r(n-i,0)$ for $0 \le i \le n-1$. But $r(n-i,0) = \Lambda_{f^{(m)},q}\left(f^{(m)}_{m(n-i)}(x,q)\right) = 0.$

As $q-$ analogue of (1.8) we get (cf. [13])

$$\sum_{n\ge 0} q^{-2\binom{n}{2}} \frac{x^n}{(-x;q^{-1})_{2n+1}} C_n(q) = 1. \tag{2.9}$$

To prove this observe that by (1.35) $\dfrac{1}{(x;q)_{n+1}} = \sum_{k\ge 0}\begin{bmatrix} n+k-1 \\ k \end{bmatrix} x^k$ and that

$$\sum_{k=0}^{n} q^{2\binom{k}{2}} \begin{bmatrix} 2n-k \\ k \end{bmatrix} (-1)^k C_{n-k}(q) = \Lambda_{f^{(2)},q}\left(f^{(2)}_{2n}(x,q)\right) = [n=0].$$

Therefore we get

$$\sum_{j\ge 0} q^{-2\binom{j}{2}} \frac{x^j}{(-x;q^{-1})_{2j+1}} C_j(q) = \sum_{j\ge 0} q^{-2\binom{j}{2}} C_j(q) x^j \sum_{k\ge 0} \begin{bmatrix} 2j+k \\ k \end{bmatrix}_{q^{-1}} (-x)^k$$

$$= \sum_n x^n \sum_{j+k=n} \begin{bmatrix} 2j+k \\ k \end{bmatrix}_q q^{2k^2-2kn} C_{n-k}(q) q^{-2\binom{n-k}{2}} (-1)^k = \sum_n x^n q^{-2\binom{n}{2}} \sum_{k=0}^{n} q^{2\binom{k}{2}} \begin{bmatrix} 2n-k \\ k \end{bmatrix} (-1)^k C_{n-k}(q) = 1.$$

In the same way we see that

$$\sum_{n\ge 0} q^{-m\binom{n}{2}} \frac{x^n}{(-x;q^{-1})_{mn+1}} C_n^{(m)}(q) = \sum_{j\ge 0} q^{-m\binom{j}{2}} C_n^{(m)}(q) x^j \sum_{k\ge 0} \begin{bmatrix} mj+k \\ k \end{bmatrix}_{q^{-1}} (-x)^k$$

$$= \sum_n x^n \sum_{j+k=n} \begin{bmatrix} mj+k \\ k \end{bmatrix}_q q^{-mjk} C_{n-k}^{(m)}(q) q^{-m\binom{n-k}{2}} (-1)^k = \sum_n x^n q^{-m\binom{n}{2}} \sum_{k=0}^{n} q^{m\binom{k}{2}} \begin{bmatrix} 2n-k \\ k \end{bmatrix} (-1)^k C_{n-k}^{(m)}(q) = 1.$$

and therefore

$$\sum_{n\ge 0} q^{-m\binom{n}{2}} \frac{x^n}{(-x;q^{-1})_{mn+1}} C_n^{(m)}(q) = 1. \tag{2.10}$$

It seems that there are no $q-$ analogues of $l_n^{(m)}(x)$ with simple recurrence relations. But there is a rather curious class of polynomials which satisfies an operator recurrence relation.

### 3. Some curious q-analogues

**3.1.** Let us consider the polynomials (cf. [8])

$$F_n(x,q) = \sum_{k=0}^{\lfloor \frac{n}{2} \rfloor} q^{\binom{k+1}{2}} \begin{bmatrix} n-k \\ k \end{bmatrix} (-1)^k x^{n-2k}. \tag{3.1}$$

They are not only $q-$ analogues of $f_n(x)$ in the sense that $\lim_{q \to 1} F_n(x,q) = f_n(x)$ but can be obtained from $f_n(x)$ by first computing the operator $f_n(x+(1-q)D_q)$ and then applying it to the constant polynomial 1.
Thus

$$F_n(x,q) = f_n(x+(1-q)D_q)1. \tag{3.2}$$

They satisfy the recurrence relation

$$F_n(x,q) = \left(x+(1-q)D_q\right)F_{n-1}(x,q) - F_{n-2}(x,q) \tag{3.3}$$

with initial values $F_0(x,q) = 1$ and $F_1(x,q) = x.$

They also satisfy

$$F_n(x,q) = xF_{n-1}(x,q) - q^{n-1}xF_{n-3}(x,q) + q^{n-1}F_{n-4}(x,q). \tag{3.4}$$

Let $\Lambda_{F,q}\left(F_n(x,q)\right) = [n=0].$

The polynomials $F_n(x,q)$ are not orthogonal. For example $\Lambda_{F,q}\left(xF_3(x,q)\right) = (q-1)q^3 \neq 0.$

Nevertheless there a very nice $q-$ analogue of (1.4):

$$x^n = \sum_{k=0}^{\lfloor \frac{n}{2} \rfloor} \left(\begin{bmatrix} n \\ k \end{bmatrix} - \begin{bmatrix} n \\ k-1 \end{bmatrix}\right)_{n-2k} F_{n-2k}(x,q). \tag{3.5}$$

This implies that the moments $\Lambda_{F,q}\left(x^{2n}\right)$ are

$$\Lambda_{F,q}\left(x^{2n}\right) = \left(\begin{bmatrix} 2n \\ n \end{bmatrix} - \begin{bmatrix} 2n \\ n-1 \end{bmatrix}\right) = q^n c_q(n), \tag{3.6}$$

where $c_q(n) = \frac{1}{[n+1]}\begin{bmatrix} 2n \\ n \end{bmatrix}$ is a explicit $q-$ analogue of the Catalan numbers.

I do not know a simple $q-$ analogue of the generating function (1.7), but we have instead

$$\sum_{n\geq 0} c_n(q) q^{-\binom{n}{2}} \frac{u^n}{\left(-q^{-n}u;q\right)_{2n+1}} = 1 \tag{3.7}$$

which is a $q-$ analogue of (1.8).

**3.2.** Let now

$$l_n(x,q) = F_n(x,q) - F_{n-2}(x,q) \tag{3.8}$$

for $n \geq 2$ and $l_0(x,q) = 1$ and $l_1(x,q) = x$.

Then

$$l_n(x,q) = \sum_{k=0}^{n} (-1)^k q^{\binom{k}{2}} \frac{[n]}{[n-k]} \begin{bmatrix} n-k \\ k \end{bmatrix} x^{n-2k} \tag{3.9}$$

for $n > 0$ and $l_0(x,q) = 1$.

The polynomials $l_n(x,q)$ satisfy

$$l_n(x,q) = \left(x + (1-q)D_q\right) l_{n-1}(x,q) - \lambda_{n-2} l_{n-2}(x,q) \tag{3.10}$$

with initial values $l_0(x,q) = 1$ and $l_1(x,q) = x$. Here $\lambda_0 = 2$ and $\lambda_n = 1$ for $n > 0$.

This can also be written as

$$l_n(x,q) = l_n(x + (1-q)D_q)1. \tag{3.11}$$

The polynomials $l_n(x,q)$ are not orthogonal.

The identity

$$x^n = \sum_{k=0}^{\left\lfloor \frac{n}{2} \right\rfloor} \begin{bmatrix} n \\ k \end{bmatrix} l_{n-2k}(x,q) \tag{3.12}$$

implies that

$$\Lambda_{l,q}\left(x^{2n}\right) = \begin{bmatrix} 2n \\ n \end{bmatrix}, \tag{3.13}$$

if we define the linear functional $\Lambda_{l,q}$ by $\Lambda_{l,q}(l_n(x,q)) = [n=0]$.

I do not know a simple $q-$ analogue of (1.16) for $b_q(u) = \sum_{n\geq 0} \begin{bmatrix} 2n \\ n \end{bmatrix} u^n$.

Instead of this

$$\sum_{n\geq 0} \begin{bmatrix} 2n \\ n \end{bmatrix} q^{-\binom{n+1}{2}} \frac{u^n}{\left(-q^{-n}u;q\right)_{2n+1}} = \sum_{n\geq 0} q^{-\binom{n+1}{2}} u^n \tag{3.14}$$

is a $q-$ analogue of $\sum_{n\geq 0} \binom{2n}{n} \frac{u^n}{(1+u)^{2n+1}} = \frac{1}{1-u}.$

### 3.3. Proofs and remarks

The polynomials $F_n(x,q)$ and $l_n(x,q)$ have been systematically studied in [8]. To prove (3.1) it suffices to compare coefficients in (3.3).

Since these polynomials are not orthogonal and thus do not satisfy a $3-$ term recurrence of the form (1.38) the above combinatorial interpretation fails.

But formula (3.3) implies Binet-type formulae for these polynomials:

Let $A$ be the operator $A = x + (1-q)D_q$. For each polynomial $p(x)$ in $x$ we define $\Phi(p(x)) = p(A)1$.

Thus $F_n(x,q) = \Phi(f_n(x)) = \Phi\left(\dfrac{\alpha^{n+1} - \beta^{n+1}}{\alpha - \beta}\right)$ and analogously

$l_n(x,q) = \Phi(l_n(x)) = \Phi(\alpha^n + \beta^n)$

for $n > 0$.

This is an exact version of a symbolic method which I used in [10].

This implies

$$\Phi\left(\sum_{k=0}^{n} \begin{bmatrix} n \\ k \end{bmatrix} \alpha^k \beta^{n-k}\right) = x^n. \tag{3.15}$$

From (1.30) we get $\sum_{k=0}^{n} \begin{bmatrix} n \\ k \end{bmatrix} \alpha^k \beta^{n-k} = x \sum_{k} \begin{bmatrix} n-1 \\ k \end{bmatrix} \alpha^k \beta^{n-k} + (q^{n-1} - 1) \sum_{k} \begin{bmatrix} n-2 \\ k \end{bmatrix} \alpha^k \beta^{n-k}.$

Since these are by induction polynomials in $x$ we get again by induction

$$\Phi\left(\sum_{k=0}^{n} \begin{bmatrix} n \\ k \end{bmatrix} \alpha^k \beta^{n-k}\right) = Ax^{n-1} + (q^{n-1} - 1)x^{n-2} = x^n.$$

To prove (3.12) observe that for odd $n$

$$\sum_{k=0}^{\lfloor n/2 \rfloor} \begin{bmatrix} n \\ k \end{bmatrix} l_{n-2k}(x,q) = \Phi\left(\sum_{k=0}^{\lfloor n/2 \rfloor} \begin{bmatrix} n \\ k \end{bmatrix} (\alpha^{n-2k} + \beta^{n-2k})\right) = \Phi\left(\sum_{k=0}^{\lfloor n/2 \rfloor} \begin{bmatrix} n \\ k \end{bmatrix} (\alpha^{n-k}\beta^k + \beta^{n-k}\alpha^k)\right)$$

$$= \Phi\left(\sum_{k=0}^{n} \begin{bmatrix} n \\ k \end{bmatrix} \alpha^k \beta^{n-k}\right) = x^n.$$

If $n = 2m$ then

$$\sum_{k=0}^{m} \begin{bmatrix} 2m \\ k \end{bmatrix} l_{2m-2k}(x,q) = \begin{bmatrix} 2m \\ m \end{bmatrix} + \Phi\left(\sum_{k=0}^{m-1} \begin{bmatrix} 2m \\ k \end{bmatrix} (\alpha^{2m-2k} + \beta^{2m-2k})\right) = \Phi\left(\sum_{k=0}^{2m} \begin{bmatrix} 2m \\ k \end{bmatrix} \alpha^k \beta^{2m-k}\right) = x^{2m}.$$

Since $l_n(x,q) = F_n(x,q) - F_{n-2}(x,q)$ we also get (3.5).

Comparing coefficients we see that (3.7) is equivalent with

$$q^{-\binom{n}{2}}\sum_{j=0}^{n}(-1)^j c_{n-j}(q)q^{\binom{j}{2}}\begin{bmatrix}2n-j\\j\end{bmatrix}=[n=0].$$

But this is clear since

$$q^n\sum_{j=0}^{n}(-1)^j c_{n-j}(q)q^{\binom{j}{2}}\begin{bmatrix}2n-j\\j\end{bmatrix}=\Lambda_{F,q}\left(\sum_{j=0}^{n}(-1)^j q^{\binom{j+1}{2}}\begin{bmatrix}2n-j\\j\end{bmatrix}x^{2n-2j}\right)=\Lambda_{F,q}\left(F_{2n}(x,q)\right)=[n=0].$$

In the same way we prove identity (3.14):

$$\sum_{n\geq 0}\begin{bmatrix}2n\\n\end{bmatrix}q^{-\binom{n+1}{2}}\frac{u^n}{\left(-q^{-n}u;q\right)_{2n+1}}=\sum_{j\geq 0}\begin{bmatrix}2j\\j\end{bmatrix}q^{-\binom{j+1}{2}}u^j\sum_{k\geq 0}\begin{bmatrix}2j+k\\k\end{bmatrix}(-1)^k q^{-jk}u^k$$

$$=\sum_{n\geq 0}u^n\sum_{k=0}^{n}(-1)^k\begin{bmatrix}2n-2k\\n-k\end{bmatrix}\begin{bmatrix}2n-k\\k\end{bmatrix}q^{-\binom{n-k+1}{2}-(n-k)k}=\sum_{n\geq 0}u^n q^{-\binom{n+1}{2}}\sum_{k=0}^{n}(-1)^k q^{\binom{k+1}{2}}\begin{bmatrix}2n-k\\k\end{bmatrix}\begin{bmatrix}2n-2k\\n-k\end{bmatrix}$$

$$=\sum_{n\geq 0}u^n q^{-\binom{n+1}{2}}\Lambda_{l,q}\left(F_{2n}(x,q)\right)$$

Since $F_n(x,q)=l_n(x,q)+F_{n-2}(x,q)$ we get $F_n(x,q)=\sum_{k=0}^{\lfloor n/2\rfloor}l_{n-2k}(x,q)$. This implies $\Lambda_{l,q}\left(F_{2n}(x,q)\right)=1$.

**3.4.** The polynomials $F_n(x,q)$ can be generalized to

$$F_n^{(m)}(x,q)=\sum_{k=0}^{\lfloor n/m\rfloor}(-1)^k q^{\binom{k+1}{2}}\begin{bmatrix}n-(m-1)k\\k\end{bmatrix}x^{n-mk} \tag{3.16}$$

which satisfy

$$F_n^{(m)}(x,q)=xF_{n-1}^{(m)}(x,q)+(1-q)D_q F_{n-m+1}^{(m)}(x,q)-F_{n-m}^{(m)}(x,q).$$

Note that the operator $(1-q)D_q$ is applied to $F_{n-m+1}^{(m)}(x,q)$.

This follows since the coefficient of $(-1)^k q^{\binom{k+1}{2}}x^{n-mk}$ of the right-hand side is

$$\begin{bmatrix} n-1-(m-1)k \\ k \end{bmatrix} - q^{-k} \begin{bmatrix} n-(m-1)k \\ k-1 \end{bmatrix} (1-q^{n-mk+1}) + q^{-k} \begin{bmatrix} n-(m-1)k-1 \\ k-1 \end{bmatrix}$$

$$= q^{-k} \left( q^k \begin{bmatrix} n-1-(m-1)k \\ k \end{bmatrix} + \begin{bmatrix} n-1-(m-1)k \\ k \end{bmatrix} \right) - q^{-k}(1-q^k) \begin{bmatrix} n-(m-1)k \\ k \end{bmatrix}$$

$$= q^{-k} \begin{bmatrix} n-(m-1)k \\ k \end{bmatrix} - q^{-k}(1-q^k) \begin{bmatrix} n-(m-1)k \\ k \end{bmatrix} = \begin{bmatrix} n-(m-1)k \\ k \end{bmatrix}$$

For $m=1$ this reduces to

$$F_n^{(1)}(x,q) = \sum_{k=0}^{n} (-1)^k q^{\binom{k+1}{2}} \begin{bmatrix} n \\ k \end{bmatrix} x^{n-k} = (x-q)(x-q^2) \cdots (x-q^n).$$

Let

$$c_n^{(m)}(q) = \Lambda_{F^{(m)},q}(x^{mn}). \tag{3.17}$$

Note that $c_n^{(2)}(q) = q^n c_n(q)$.

If we apply $\Lambda_{F^{(m)},q}$ to $F_{mn}^{(m)}(x,q)$ we get a recurrence for $c_n^{(m)}(q)$.

$$\sum_{k=0}^{n} (-1)^k q^{\binom{k+1}{2}} \begin{bmatrix} mn-(m-1)k \\ k \end{bmatrix} c_{n-k}^{(m)}(q) = 0. \tag{3.18}$$

The numbers $c_n^{(m)}(q)$ satisfy

$$\sum_{n \geq 0} c_n^{(m)}(q) q^{-\binom{n+1}{2}} \frac{x^n}{(-q^{-n}x;q)_{mn+1}} = 1. \tag{3.19}$$

For the left-hand side is

$$\sum_{j \geq 0} c_j^{(m)}(q) q^{-\binom{j+1}{2}} x^j \sum_{k \geq 0} \begin{bmatrix} mj+k \\ k \end{bmatrix} (-x)^k q^{-kj} = \sum_{n \geq 0} x^n \sum_{k=0}^{n} (-1)^k c_{n-k}^{(m)}(q) q^{-\binom{n-k+1}{2}-k(n-k)} \begin{bmatrix} mn-(m-1)k \\ k \end{bmatrix}$$

$$= \sum_{n \geq 0} x^n q^{-\binom{n+1}{2}} \sum_{k=0}^{n} (-1)^k q^{\binom{k+1}{2}} \begin{bmatrix} mn-(m-1)k \\ k \end{bmatrix} c_{n-k}^{(m)}(q) = \sum_{n \geq 0} x^n q^{-\binom{n+1}{2}} \Lambda_{F^{(m)},q}\left(F_{mn}^{(m)}(x,q)\right) = 1.$$

In the same way the polynomials $l_n(x,q) = L_n^{(2)}(x,q)$ can be generalized to

$$L_n^{(m)}(x,q) = \sum_{k=0}^{\lfloor n/m \rfloor} (-1)^k q^{\binom{k}{2}} \begin{bmatrix} n-(m-1)k \\ k \end{bmatrix} \frac{[n-(m-2)k]}{[n-(m-1)k]} x^{n-mk}$$

which satisfy
$$L_n^{(m)}(x,q) = xL_{n-1}^{(m)}(x,q) + (1-q)D_q L_{n-m+1}^{(m)}(x,q) - \lambda_{n-m}L_{n-m}^{(m)}(x,q) \tag{3.20}$$

with $\lambda_0 = 2$ and $\lambda_n = 1$.

As above we have
$$L_n^{(m)}(x,q) = F_n^{(m)}(x,q) - F_{n-m}^{(m)}(x,q). \tag{3.21}$$

Let
$$b_n^{(m)}(q) = \Lambda_{L^{(m)},q}\left(x^{mn}\right). \tag{3.22}$$

Then we get the recurrence
$$\sum_{k=0}^{n}(-1)^k q^{\binom{k}{2}}\begin{bmatrix}mn-(m-1)k\\k\end{bmatrix}\frac{[mn-(m-2)k]}{[mn-(m-1)k]}b_{n-k}^{(m)}(q) = 0. \tag{3.23}$$

Generalizing (3.14) we get
$$\sum_{n\geq 0}b_n^{(m)}(q)q^{-\binom{n+1}{2}}\frac{u^n}{\left(-q^{-n}u;q\right)_{mn+1}} = \sum_{n\geq 0}u^n q^{-\binom{n+1}{2}}. \tag{3.24}$$

Since by (3.21) we have $F_n^{(m)}(x,q) = \sum_{k=0}^{\left\lfloor\frac{n}{m}\right\rfloor}L_{n-mk}^{(m)}(x,q)$ we have $\Lambda_{L^{(m)},q}\left(F_n^{(m)}(x,q)\right) = 1$.

Thus we get
$$\sum_{n\geq 0}b_n^{(m)}(q)q^{-\binom{n+1}{2}}\frac{u^n}{\left(-q^{-n}u;q\right)_{mn+1}} = \sum_{j\geq 0}b_j^{(m)}(q)q^{-\binom{j+1}{2}}u^j\sum_{k\geq 0}\begin{bmatrix}mj+k\\k\end{bmatrix}(-1)^k q^{-jk}u^k$$
$$= \sum_{n\geq 0}u^n\sum_{k=0}^{n}(-1)^k b_{n-k}^{(m)}(q)\begin{bmatrix}mn-(m-1)k\\k\end{bmatrix}q^{-\binom{n-k+1}{2}-(n-k)k}$$
$$= \sum_{n\geq 0}u^n q^{-\binom{n+1}{2}}\sum_{k=0}^{n}(-1)^k q^{\binom{k+1}{2}}\begin{bmatrix}mn-(m-1)k\\k\end{bmatrix}b_{n-k}^{(m)}(q) = \sum_{n\geq 0}u^n q^{-\binom{n+1}{2}}\Lambda_{L^{(m)},q}\left(F_{n-k}^{(m)}(x,q)\right) = \sum_{n\geq 0}u^n q^{-\binom{n+1}{2}}.$$

For $m=1$ we get $L_n^{(1)}(x,q) = (x-1-q^n)\prod_{j=1}^{n-1}\left(x-q^j\right)$

for $n \geq 1$.

## 4. q-Chebyshev polynomials

Now we come to a class of orthogonal polynomials where almost all facts from the classical case have simple counterparts.

**4.1.** The polynomials

$$u_n(x,q) = \sum_{k=0}^{\lfloor \frac{n}{2} \rfloor} \begin{bmatrix} n-k \\ k \end{bmatrix} q^{k^2} (-1)^k \frac{x^{n-2k}}{(-q;q)_k (-q^{n+1-k};q)_k} \qquad (4.1)$$

will be called special $q-$ Chebyshev polynomials of the second kind.

They satisfy the recurrence relation

$$u_n(x,q) = x u_{n-1}(x,q) - \frac{q^{n-1}}{(1+q^{n-1})(1+q^n)} u_{n-2}(x,q) \qquad (4.2)$$

with initial values $u_0(x,q) = 1$ and $u_1(x,q) = x$.

The polynomials $u_n(x,q)$ are orthogonal with respect to the linear functional defined by

$$\Lambda_{u,q}(u_n(x,q)) = [n=0]. \qquad (4.3)$$

More precisely we have

$$\Lambda_{u,q}(u_n(x,q) u_m(x,q)) = \frac{q^{\binom{n+1}{2}}}{(-q;q)_n (-q^2;q)_n} [n=m]. \qquad (4.4)$$

The identity

$$\sum_{k=0}^{\lfloor \frac{n}{2} \rfloor} \frac{\begin{bmatrix} n \\ k \end{bmatrix} - \begin{bmatrix} n \\ k-1 \end{bmatrix}}{(-q;q)_k (-q^{n+2-2k};q)_k} u_{n-2k}(x,q) = x^n \qquad (4.5)$$

gives the moments

$$\Lambda_{u,q}(x^{2n}) = \frac{1}{[n+1]} \begin{bmatrix} 2n \\ n \end{bmatrix} \frac{q^n}{(-q;q)_n (-q^2;q)_n} = (-1)^n q^{n^2}(1+q) \begin{bmatrix} \frac{1}{2} \\ n+1 \end{bmatrix}_{q^2} = \mathbf{C}_n(q) \qquad (4.6)$$

where $\mathbf{C}_n(q)$ is a $q-$ Catalan number in the sense of Andrews [2].

As $q-$ analogue of $C\left(\dfrac{u}{4}\right)=1+\dfrac{u}{4}C\left(\dfrac{u}{4}\right)^2$ we get for the generating function
$\mathbf{C}(u,q)=\sum_{n\geq 0}\mathbf{C}_n(q)u^n$

$$\frac{\mathbf{C}(u,q)+q\mathbf{C}(qu,q)}{1+q}=1+\frac{qu}{(1+q)^2}\mathbf{C}(u,q)\mathbf{C}(qu,q). \tag{4.7}$$

Let $h(u):=\dfrac{(u;q^2)_\infty}{(qu;q^2)_\infty}$. This is a $q-$ analogue of $\sqrt{1-u}$ since $h(u)h(qu)=\dfrac{(u;q^2)_\infty}{(q^2u;q^2)_\infty}=1-u$.

The formula

$$\mathbf{C}(u,q)=(1+q)\frac{1-h(u)}{u} \tag{4.8}$$

is a $q-$ analogue of $\sum_{n\geq 0}\dfrac{C_n}{4^n}u^n=2\dfrac{1-\sqrt{1-u}}{u}$.

**4.2.** The special $q-$ Chebyshev polynomials of the first kind are the polynomials

$$t_n(x,q)=\sum_{k=0}^{\left\lfloor\frac{n}{2}\right\rfloor}(-1)^k q^{k^2}\frac{[n]}{[n-k]}\begin{bmatrix}n-k\\k\end{bmatrix}\frac{1}{(-q;q)_k(-q^{n-k};q)_k}x^{n-2k}. \tag{4.9}$$

The polynomials $t_n(x,q)$ satisfy the recurrence relation

$$t_n(x,q)=xt_{n-1}(x,q)-\lambda_{n-2}(q)t_{n-2}(x,q) \tag{4.10}$$

with $\lambda_0(q)=\dfrac{q}{1+q}$ and $\lambda_n(q)=\dfrac{q^{n+1}}{(1+q^n)(1+q^{n+1})}$ for $n>0$.

It is easy to verify that $t_0(x,q)=u_0(x,q)=1$, $t_1(x,q)=u_1(x,q)=x$ and for $n\geq 2$

$$t_n(x,q)=u_n(x,q)-\frac{q^{2n-1}}{(1+q^{n-1})(1+q^n)}u_{n-2}(x,q). \tag{4.11}$$

The polynomials $t_n(x,q)$ are orthogonal with respect to the linear functional $\Lambda_{t,q}$ defined by

$$\Lambda_{t,q}(t_n(x,q))=[n=0]. \tag{4.12}$$

More precisely we have for $n>0$

$$\Lambda_{t,q}(t_n(x,q)t_m(x,q))=\frac{q^{\binom{n+1}{2}}}{(-q;q)_{n-1}(-q;q)_n}[n=m]. \tag{4.13}$$

The identity ([11], Theorem 4.3)

$$x^n = \sum_{k=0}^{\lfloor \frac{n}{2} \rfloor} \begin{bmatrix} n \\ k \end{bmatrix} \frac{q^k}{(-q;q)_k \left(-q^{n-2k+1};q\right)_k} t_{n-2k}(x,q) \tag{4.14}$$

implies the moments

$$\Lambda_{t,q}\left(x^{2n}\right) = \begin{bmatrix} 2n \\ n \end{bmatrix} \frac{q^n}{(-q;q)_n^2}. \tag{4.15}$$

Let

$$G(u,q) := \sum_{k \geq 0} \frac{q^{k^2-k}}{(q^2;q^2)_k}(-1)^k u^k = \left(u;q^2\right)_\infty. \tag{4.16}$$

By the $q-$ binomial theorem

$$g(u) = \frac{G(qu,q)}{G(u,q)} = \sum_{n \geq 0} \begin{bmatrix} 2n \\ n \end{bmatrix} \frac{1}{(-q;q)_n^2} u^n. \tag{4.17}$$

Note that $g(u)$ is a $q-$ analogue of $\frac{1}{\sqrt{1-u}}$ since $G(u,q) = \left(u;q^2\right)_\infty$ implies $g(u)g(qu) = \frac{1}{1-u}$. The generating function of the moments is

$$\sum_{n \geq 0} \Lambda_{t,q}\left(x^{2n}\right) u^n = g(qu) = \frac{G(q^2 u, q)}{G(qu, q)}. \tag{4.18}$$

**Notes**

In [11] we introduced bivariate $q-$Chebyshev polynomials $T_n(x,s,q)$ of the first kind
by $T_n(x,s,q) = \left(1+q^{n-1}\right)xT_{n-1}(x,s,q) + q^{n-1}sT_{n-2}(x,s,q)$
with initial values $T_0(x,s,q) = 1$ and $T_1(x,s,q) = x$
and bivariate $q-$Chebyshev polynomials $U_n(x,s,q)$ of the second kind
by $U_n(x,s,q) = \left(1+q^n\right)xU_{n-1}(x,s,q) + q^{n-1}sU_{n-2}(x,s,q)$
with initial values $U_0(x,s,q) = 1$ and $U_1(x,s,q) = (1+q)x.$

We then have

$$t_n(x,q) = \frac{T_n(x,-1,q)}{(-q;q)_{n-1}} \tag{4.19}$$

for $n > 0$ and

$$u_n(x,q) = \frac{U_n(x,-1,q)}{(-q;q)_n}. \tag{4.20}$$

Similar polynomials have also appeared in other publications, cf. [11] or [12] and the literature cited there. They are related to the Al-Salam and Ismail polynomials introduced in [1].

The recurrence relations can be easily verified by comparing coefficients. Proofs for (4.5) and (4.14) can be found in [11].

Formula (4.8) follows from the $q-$ binomial theorem (1.31) since

$$h(u) = \frac{(u;q^2)_\infty}{(qu;q^2)_\infty} = \sum_{n\geq 0} \frac{(q^{-1};q^2)_n}{(q^2;q^2)_n}(qu)^n = 1 + \sum_{n\geq 0} q^{n+1}\frac{(q^{-1};q^2)_{n+1}}{(q^2;q^2)_{n+1}}u^{n+1}$$

$$= 1 - \frac{u}{1+q}\sum_{n\geq 0}\frac{1}{[n+1]}\begin{bmatrix}2n\\n\end{bmatrix}\frac{q^n}{(-q;q)_n(-q^2;q)_n}u^n = 1 - \frac{u}{1+q}\sum_{n\geq 0}\mathbf{C}_n(q)u^n.$$

(4.8) implies $\dfrac{u^2 q \mathbf{C}(u,q)\mathbf{C}(qu,q)}{(1+q)^2} = (1-h(u))(1-h(qu)) = \dfrac{u}{1+q}(\mathbf{C}(u,q)+q\mathbf{C}(qu,q))-u$

and thus (4.7).

The $q-$ binomial theorem gives

$$g(qu) = \frac{(q^2u;q^2)_\infty}{(qu;q^2)_\infty} = \sum_{k\geq 0}\frac{(q;q^2)_k}{(q^2;q^2)_k}(qu)^k = \sum_{k\geq 0}\begin{bmatrix}2k\\k\end{bmatrix}\frac{q^k}{(-q;q)_k^2}u^k.$$

Note that the binomial theorem gives $\dfrac{1}{G(u,q)} = \dfrac{1}{(u;q^2)_\infty} = \sum_{k\geq 0}\dfrac{u^k}{(q^2;q^2)_k}.$

Thus by comparison of coefficients (4.17) is equivalent with the well-known formula

$$\sum_{j=0}^n q^{j^2}\begin{bmatrix}n\\j\end{bmatrix}^2 = \begin{bmatrix}2n\\n\end{bmatrix}. \tag{4.21}$$

## 5. A slight extension.

**5.1.** Consider the orthogonal polynomials (1.38) with $\lambda_k = \lambda_k(z,q) = \dfrac{q^{k+1}}{(1-q^k z)(1-q^{k+1}z)}.$

Calling them $f_n(x,z,q)$ we get

$$f_n(x,z,q) = \sum_{k=0}^{\left\lfloor \frac{n}{2} \right\rfloor} q^{k^2} \begin{bmatrix} n-k \\ k \end{bmatrix} \frac{(-1)^k}{(z;q)_k \left(q^{n-k}z;q\right)_k} x^{n-2k}. \tag{5.1}$$

Note that $f_n(x,-q,q) = u_n(x,q).$

These polynomials are also related to the Al-Salam and Ismail polynomials (cf. [1],[11] or [12] and the literature cited there).

By (1.45) we get for the generating functions $\Phi_f(u,z,q) = \sum_n c(2n,0,z,q)u^n$ and $\Psi_f(u,z,q) = \sum_n b(2n,0,z,q)u^n$

$$\Phi_f(u,z,q) = 1 + \frac{qu}{(1-z)(1-qz)} \Phi_f(u,z,q)\Psi_f(u,z,q).$$

Since $\lambda_{k+1}(z,q) = \frac{q^{k+2}}{(1-q^{k+1}z)(1-q^{k+2}z)} = q\lambda_k(qz,q)$ we have $\Psi_f(u,z,q) = \Phi_f(qu,qz,q).$

Therefore $\Phi_f(u,z,q)$ satisfies

$$\Phi_f(u,z,q) = 1 + \frac{qu}{(1-z)(1-qz)} \Phi_f(u,z,q)\Phi_f(qu,qz,q). \tag{5.2}$$

For $q \to 1$ this gives $\Phi_f(u,z,1) = 1 + \frac{qu}{(1-z)^2} \Phi_f(u,z,1)^2$ and thus

$$\Phi_f(u,z,1) = C\left(\frac{u}{(1-z)^2}\right). \tag{5.3}$$

In the general case there are no simple formulae for $c(2n,0,z,q),$ but there is a simple representation for their generating functions.

Let

$$G(u,z,q) = \sum_{n \geq 0} q^{n^2-n} \frac{u^n}{(z;q)_n (q;q)_n} \tag{5.4}$$

which is a $q-$ analogue of the exponential series $\sum_{n \geq 0} \frac{1}{n!}\left(\frac{u}{1-z}\right)^n.$

Then

$$\psi(u,z,q) = \frac{G(qu,qz,q)}{G(u,z,q)} \tag{5.5}$$

satisfies

$$\psi(u,z,q) = 1 - \frac{u}{(1-z)(1-qz)}\psi(u,z,q)\psi(qu,qz,q). \tag{5.6}$$

Therefore we get

$$\Phi_f(u,z,q) = \psi(-qu,z,q) = \frac{G(-q^2u,qz,q)}{G(-qu,z,q)}. \tag{5.7}$$

This follows from

$$G(qu,qz,q) - G(u,z,q) = \sum_n q^{n^2-n}\left(\frac{q^n u^n}{(qz;q)_n (q;q)_n} - \frac{u^n}{(z;q)_n (q;q)_n}\right)$$

$$= \sum_n q^{n^2-n}\frac{u^n\left(q^n(1-z) - (1-q^n z)\right)}{(z;q)_{n+1}(q;q)_n} = -\frac{u}{(1-z)(1-qz)}\sum_n q^{n^2-n}\frac{q^{2n}u^n}{(q^2 z;q)_n (q;q)_n}$$

$$= -\frac{u}{(1-z)(1-qz)} G(q^2 u, q^2 z, q)$$

which implies

$$\psi(u,z,q) - 1 = \frac{G(qu,qz,q)}{G(u,z,q)} - 1 = -\frac{u}{(1-z)(1-qz)}\frac{G(q^2u,q^2z,q)}{G(u,z,q)}$$

$$= -\frac{u}{(1-z)(1-qz)}\frac{G(q^2u,q^2z,q)}{G(qu,qz,q)}\frac{G(qu,qz,q)}{G(u,z,q)} = -\frac{u}{(1-z)(1-qz)}\psi(u,z,q)\psi(qu,qz,q).$$

It is clear that $\tilde{\psi}(u,z,q) = \psi(uz,z,q)$ which satisfies

$$\tilde{\psi}(u,z,q) = 1 - \frac{uz}{(1-z)(1-qz)}\tilde{\psi}(u,z,q)\tilde{\psi}(u,qz,q) \tag{5.8}$$

is equivalent with

$$\tilde{\psi}(u,z,q) = \frac{G(-quz,qz,q)}{G(-uz,z,q)}. \tag{5.9}$$

**5.2.** If we choose

$$\lambda_0(z,q) = \frac{q}{1-z}$$

$$\lambda_n(z,q) = \frac{q^{n+1}}{(1-q^{n-1}z)(1-q^n z)} \tag{5.10}$$

we get the polynomials

$$l_n(x,z,q) = \sum_{k=0}^{\lfloor\frac{n}{2}\rfloor} (-1)^k q^{k^2} \frac{\left(\begin{bmatrix}n-k\\k\end{bmatrix} - q^{n-k-1}z\begin{bmatrix}n-k-1\\k-1\end{bmatrix}\right)}{(z;q)_k (q^{n-1-k}z;q)_k} x^{n-2k}. \tag{5.11}$$

Note that $t_n(x,q) = l_n(x,-q,q)$.

For the generating function we get

$$\Phi_l(u,z,q) = 1 + \frac{qu}{1-z}\Phi_l(u,z,q)\Psi_l(u,z,q) \tag{5.12}$$

and

$$\Psi_l(u,z,q) = 1 + \frac{q^2 u}{(1-z)(1-qz)}\Psi_l(u,z,q)\Psi_l(qu,qz,q). \tag{5.13}$$

This implies

$$\Psi_l(u,z,q) = \Phi_f(qu,z,q^2,q) = \frac{G(-q^3 u, qz, q)}{G(-q^2 u, z, q)}. \tag{5.14}$$

Since

$$G(qu,z,q) - G(u,z,q) = \sum_k q^{k^2-k} \frac{u^k(q^k-1)}{(z;q)_k (q;q)_k}$$

$$= -\frac{u}{1-z}\sum_k q^{k^2-k}\frac{(q^2 u)^k}{(qz;q)_k (q;q)_k} = -\frac{u}{1-z} G(q^2 u, qz, q)$$

we get for $\chi(u,z,q) = \dfrac{G(-q^2 u, z, q)}{G(-qu, z, q)}$

$$\chi(u,z,q) = \frac{G(-q^2 u, z, q)}{G(-qu, z, q)} = 1 + \frac{qu}{1-z}\frac{G(-q^3 u, qz, q)}{G(-qu, z, q)} = 1 + \frac{qu}{1-z}\frac{G(-q^3 u, qz, q)}{G(-q^2 u, z, q)}\frac{G(-q^2 u, z, q)}{G(-qu, z, q)}$$

$$= 1 + \frac{qu}{1-z}\Psi_l(u,z,q)\chi(u,z,q).$$

(5.12) implies $\Phi_l(u,z,q) = \chi(u,z,q)$.

Thus

$$\Phi_l(u,z,q) = \frac{G(-q^2 u, z, q)}{G(-qu, z, q)}. \tag{5.15}$$

For the special case $t_n(x,q) = l_n(x,-q,q)$ we get $\Phi_l(u,-q,q) = \dfrac{G(-q^2 u, -q, q)}{G(-qu, -q, q)}$ which is the same as (4.18).

## 5.3.

Consider the series

$$F(u,z,q) = \sum_{n \geq 0} \frac{u^n}{(z;q)_n (q;q)_n} \tag{5.16}$$

and let

$$\varphi(u,z,q) = \frac{F(u,qz,q)}{F(u,z,q)}. \tag{5.17}$$

Then

$$\varphi(u,z,q) = 1 - \frac{uz}{(1-z)(1-qz)} \varphi(u,z,q)\varphi(u,qz,q). \tag{5.18}$$

This follows from

$$\varphi(u,z,q) - 1 = \frac{F(u,qz,q) - F(u,z,q)}{F(u,z,q)} = \frac{1}{F(u,z,q)} \sum_{n \geq 0} \frac{1}{(z;q)_{n+1}(q;q)_n} u^n \left((1-z) - (1-q^n z)\right)$$

$$= \frac{-zu}{(1-z)(1-qz)} \frac{F(u,q^2z,q)}{F(u,z,q)} = \frac{-zu}{(1-z)(1-qz)} \varphi(u,z,q)\varphi(u,qz,q).$$

Comparing (5.18) and (5.8) we that

$$\frac{F(u,qz,q)}{F(u,z,q)} = \frac{G(quz,qz,q)}{G(uz,z,q)}. \tag{5.19}$$

In fact we have more precisely

$$F(u,qz,q)G(uz,z,q) = F(u,z,q)G(quz,qz,q) = \sum_{n \geq} \frac{(q^n z^2;q)_n}{(z;q)_n (qz;q)_n (q;q)_n} u^n. \tag{5.20}$$

Comparing coefficients this is equivalent with

$$\sum_{k=0}^{n} \begin{bmatrix} n \\ k \end{bmatrix} q^{k^2} z^k \frac{(z;q)_n (qz;q)_n}{(z;q)_{n-k}(qz;q)_k} = \sum_{k=0}^{n} \begin{bmatrix} n \\ k \end{bmatrix} q^{k^2-k} z^k \frac{(z;q)_n (qz;q)_n}{(qz;q)_{n-k}(z;q)_k} = (q^n z^2;q)_n.$$

This follows e.g. from the $q-$ Zeilberger algorithm. We use the Mathematica implementation of PeterPaule and Axel Riese [15]:

```
qZeil[q^(k^2) z^k/qPochhammer[z, q, n - k] qBinomial[n, k, q] qPochhammer[z, q, n]
  qPochhammer[q z, q, n]/qPochhammer[q z, q, k], {k, 0, n}, n, 1]
```

$$\text{SUM}[n] == \frac{\left(1 - q^{-2+2n} z^2\right)\left(1 - q^{-1+2n} z^2\right) \text{SUM}[-1+n]}{1 - q^{-1+n} z^2}$$

and

```
qZeil[q^(k^2 - k) z^k/qPochhammer[q z, q, n - k] qBinomial[n, k, q] qPochhammer[z, q, n]
  qPochhammer[q z, q, n]/qPochhammer[ z, q, k], {k, 0, n}, n, 1]
```

$$\text{SUM}[n] == \frac{\left(1 - q^{-2+2n} z^2\right)\left(1 - q^{-1+2n} z^2\right) \text{SUM}[-1+n]}{1 - q^{-1+n} z^2}$$

By (5.19) we can express the generating function $\Phi_f(u,z,q)$ also using $F(u,z,q)$.

$$\Phi_f(u,z,q) = \frac{F\left(-\dfrac{qu}{z}, qz, q\right)}{F\left(-\dfrac{qu}{z}, z, q\right)} = \frac{G(-q^2 u, qz, q)}{G(-qu, z, q)}. \tag{5.21}$$

For the special case $z = -q$ we get another representation of the generating function of the Andrews $q-$ Catalan numbers:

$$\mathbf{C}(u,q) = \Phi_f(u,-q,q) = \frac{F(u,-q^2,q)}{F(u,-q,q)} = \frac{G(-q^2 u,-q^2,q)}{G(-qu,-q,q)}. \tag{5.22}$$

Since by (1.32) and (1.33) $F(u,-q-q)G(-u,-q,q) = 1$ this can be written as

$$\mathbf{C}(u,q) = \sum_n \frac{u^n}{(-q^2;q)_n (q;q)_n} \sum_n \frac{q^{n^2-n}(-u)^n}{(q^2;q^2)_n} = \sum_n \frac{q^n u^n}{(-q^2;q^2)_n} \sum_n \frac{q^{n^2+n}(-u)^n}{(-q^2;q)_n (q;q)_n}$$

This is equivalent with the following two (different) expressions for $\mathbf{C}_n(q)$.

$$\begin{aligned}\frac{1+q}{(q^2;q^2)_n} \sum_{k=0}^n (-1)^{n-k} q^{2\binom{n-k}{2}} \begin{bmatrix} n \\ k \end{bmatrix}_{q^2} \frac{1}{1+q^{k+1}} &= \mathbf{C}_n(q), \\ \frac{(1+q)q^n}{(q^2;q^2)_n} \sum_{k=0}^n (-1)^k q^{k^2} \begin{bmatrix} n \\ k \end{bmatrix}_{q^2} \frac{1}{1+q^{k+1}} &= \mathbf{C}_n(q).\end{aligned} \tag{5.23}$$

**Remark**

In their paper [4] M. J. Cantero and A. Iserles prove that the rational functions $a_n(z,q)$ defined by $a_0(z,q) = 1$ and

$$\sum_{j=0}^n \frac{a_{n-j}(z,q)}{(q;q)_j (z;q)_j} = \frac{q^n}{(q;q)_n (z;q)_n} \tag{5.24}$$

for $n > 0$ satisfy

$$\lim_{q \to 1} a_n(z,q) = (-1)^n C_{n-1} \frac{z^{n-1}}{(1-z)^{2n-1}}. \tag{5.25}$$

This result also follows from our considerations.

For the equations (5.24) are equivalent with

$$\sum_{k\geq 0} a_k(z,q) u^k \sum_{\ell \geq 0} \frac{1}{(z;q)_\ell (q;q)_\ell} u^\ell = \sum_{n\geq 0} \frac{1}{(z;q)_n (q;q)_n} (qu)^n \quad (5.26)$$

and thus with

$$\sum_{k\geq 0} a_k(z,q) u^k = \frac{F(qu,z,q)}{F(u,z,q)}.$$

From

$$F(qu,z,q) - F(u,z,q) = \sum_{n\geq 0} \frac{q^n u^n}{(z;q)_n (q;q)_n} - \sum_{n\geq 0} \frac{u^n}{(z;q)_n (q;q)_n} = -\sum_{n\geq 0} \frac{u^n}{(z;q)_n (q;q)_{n-1}}$$

$$= -\frac{u}{1-z} \sum_{n\geq 0} \frac{u^n}{(qz;q)_n (q;q)_n} = -\frac{u}{1-z} F(u,qz,q)$$

follows

$$\frac{F(qu,z,q)}{F(u,z,q)} = 1 - \frac{u}{1-z} \frac{F(u,qz,q)}{F(u,z,q)} = 1 - \frac{u}{1-z} \varphi(u,z,q). \quad (5.27)$$

By (5.18) $\varphi(u,z,1) = 1 - \frac{uz}{(1-z)^2} \varphi(u,z,1)^2 = C\left(-\frac{uz}{(1-z)^2}\right)$ which implies (5.25).

**Final remarks**

Some results hold also in a slightly more general version by introducing a third parameter $s$. The polynomials $f_n^{(m)}(x,q,s)$ which satisfy $f_n^{(m)}(x,q,s) = x f_{n-1}^{(m)}(x,q,s) - q^{n-m} s f_{n-m}^{(m)}(x,q,s)$ are given by the formulae $f_n^{(m)}(x,q,s) = \sum_{k=0}^{\lfloor n/m \rfloor} q^{m\binom{k}{2}} \begin{bmatrix} n-(m-1)k \\ k \end{bmatrix} (-s)^k x^{n-mk}$. Both the recurrences and the coefficients are given by closed formulae. For $m=1$ we also have a nice product representation $f_n^{(1)}(x,q) = (x-s)(x-qs)\cdots(x-q^{n-1}s)$.

As already mentioned there is no interesting $q-$ analogue of $l_n^{(m)}(x)$ with simple recurrences whose coefficients are given by a closed formula. But if we define
$l_n^{(m)}(x,q,s) = f_n^{(m)}(x,q,s) - f_{n-m}^{(m)}(x,q,q^{m-1}s)$ then we get at least a nice formula for the polynomials

$$l_n^{(m)}(x,q,s) = \sum_{k=0}^{\lfloor n/m \rfloor} q^{m\binom{k}{2}} \frac{[n-(m-2)k]}{[n-(m-1)k]} \begin{bmatrix} n-(m-1)k \\ k \end{bmatrix} (-s)^k x^{n-mk}.$$

For $m=2$ these polynomials reduce to the Carlitz $q-$ Lucas polynomials

$$\sum_{k=0}^{\lfloor n/2 \rfloor} q^{2\binom{k}{2}} \frac{[n]}{[n-k]} \begin{bmatrix} n-k \\ k \end{bmatrix} (-s)^k x^{n-2k}.$$

The polynomials $F_n^{(m)}(x,q) = \sum_{k=0}^{\lfloor n/m \rfloor} (-1)^k q^{\binom{k+1}{2}} \begin{bmatrix} n-(m-1)k \\ k \end{bmatrix} x^{n-mk}$ and

$L_n^{(m)}(x,q) = \sum_{k=0}^{\lfloor n/m \rfloor} (-1)^k q^{\binom{k}{2}} \begin{bmatrix} n-(m-1)k \\ k \end{bmatrix} \frac{[n-(m-2)k]}{[n-(m-1)k]} x^{n-mk}$

have closed formulae for the coefficients and curious $q-$ analogues of the recurrence relations.

For $m=2$ they have moreover simple closed formulae for the moments.

The most important special cases are the monic $q-$ Chebyshev polynomials which are orthogonal, have simple recurrence relations and closed formulae for both their coefficients and for the moments. It would be interesting if there exist $m-$extensions of these formulae.

Finally we have considered the polynomials $f_n(x,z,q)$ and $l_n(x,z,q)$ which also have closed formulae for their coefficients and which for $z=-q$ reduce to special Chebyshev polynomials and for $z=0$ to the Carlitz $q-$ Fibonacci polynomials $f_n^{(2)}(x,q,q)$. Their generating functions can be expressed as $\Phi_f(u,z,q) = \dfrac{G(-q^2 u, qz, q)}{G(-qu, z, q)}$ and

$\Phi_l(u,z,q) = \dfrac{G(-q^2 u, z, q)}{G(-qu, z, q)}$ respectively if $G(u,z,q) = \sum_{n \geq 0} q^{n^2-n} \dfrac{u^n}{(z;q)_n (q;q)_n}$.